\documentclass[12pt]{article}
\usepackage{amssymb}
\usepackage{amsfonts}
\usepackage{mathrsfs}
\usepackage{graphicx}
\usepackage{amsmath}
\usepackage{epsfig}
\usepackage{verbatim}
\usepackage{url}

\usepackage[top=1in, bottom=1in, left=.75in, right=.75in]{geometry}

\newcommand{\lab}[1]{\label{#1}}                

\newcommand{\remove}[1]{}
\newcommand\eqn[1]{(\ref{#1})}

\newcommand{\be}{\begin{equation}}
\newcommand{\ee}{\end{equation}}
\newcommand{\bea}{\begin{eqnarray}}
\newcommand{\eea}{\end{eqnarray}}
\newcommand{\bean}{\begin{eqnarray*}}
\newcommand{\eean}{\end{eqnarray*}}

\newtheorem{thm}{Theorem}
\newtheorem{cor}[thm]{Corollary}

\newtheorem{lemma}[thm]{Lemma}

\newtheorem{prop}[thm]{Proposition}
\def\proof{\noindent{\bf Proof.\ }  }
\def\qed{~~\vrule height8pt width4pt depth0pt}

\def\S{{\mathscr S}}

\def\G{{\mathcal G}}
\def\D{{\mathcal D}}

\def\ex{{\bf E}}
\def\pr{{\bf P}}

\def\B{{\cal B}}

\def\N{{\mathcal N}}

\def\d{d_{\max}}

\def\eps{\epsilon}
\def\la{\lambda}

\def\ss{\smallskip}

\def\no{\noindent}

\catcode`@=11 \@addtoreset{equation}{section}

\catcode`@=12
\date{}

\title{Distributions of sparse spanning subgraphs in random graphs}
\author{Pu Gao\thanks{Research supported by the Humboldt
Foundation}\\ Max-Planck-Institut f\"{u}r Informatik\\
janegao@mpi-inf.mpg.de}
\begin{document}
\maketitle

\begin{abstract}

We describe a general approach of determining the distribution of the number of
certain types of spanning subgraphs in the random graph $\G(n,p)$. Using this approach, we reprove the distribution of the number of
Hamilton cycles, with a proof that is much shorter than previously known proofs. We also achieve new results on
determining the distribution of the number of spanning triangle-free subgraphs and the number of triangle-factors.

\end{abstract}

\section{Introduction}

The distributions of subgraphs with fixed sizes in various random
graph models have been investigated by many authors. A general
approach by Ruci\'{n}ski~\cite{R,R3} showed that the numbers of
subgraphs with fixed sizes in the binomial model $\G(n,p)$ are
asymptotically normal for a large range of $p$. On the other hand,
studies of distributions of subgraphs of sizes growing with $n$, for
example, the spanning subgraphs, are much less common. The first
breakthrough is perhaps due to Robinson and Wormald~\cite{RW4,RW5}
on proving that random regular graphs are a.a.s.\ Hamiltonian. Based
on their work, Janson~\cite{J} deduced the limiting distribution of
the number of Hamilton cycles in random regular graphs. The
distributions of some types of spanning subgraphs (perfect
matchings, Hamilton cycles, spanning trees) in random  graphs
 $\G(n,p)$ and $\G(n,m)$  were determined by Janson~\cite{J3}. These
distributions behave significantly differently in $\G(n,m)$ and
$\G(n,p)$. It was shown that within a big range of $m$, the numbers
of these spanning subgraphs are asymptotically normally distributed
in $\G(n,m)$, whereas in the corresponding $\G(n,p)$ with
$p=m/\binom{n}{2}$, these random variables are asymptotically
log-normally distributed. This is because the expectations of these
variables in $\G(n,m)$ grow very fast as $m$ grows. Therefore, even
though the number of edges in $\G(n,p)$ has small deviation, the
deviation of these random variables (e.g.\ the number of perfect
matchings) can eventually be very large. This same phenomena was
observed by the author~\cite{G4} while studying the distribution of
the number of $d$-factors in $\G(n,p)$.

In this paper, we take the technique that was used in~\cite{G4} (for the study of $d$-factors) and
extend and generalise it into a method for
studying a broader class of large (spanning) subgraphs. In
Section~\ref{sec:general}, we describe the general method
(Theorems~\ref{t:Gnm} and~\ref{t:Gnp}) and give conditions under
which the distribution of the random variable under investigation
will follow a pattern of concentration in $\G(n,m)$ and log-normal
distribution in $\G(n,p)$, which we call the {\em log-normal
paradigm} in this paper. The method is also extended to cope with
probability spaces of random directed graphs (See
Theorem~\ref{t:digraph}). To show the power of the method, we reprove the
distribution of the number of Hamilton cycles. The problem on
the number of Hamilton cycles has been studied in the past by a few
authors. The first investigation was done by Wright for the directed
Hamilton cycles in~\cite{W7} and then the undirected Hamilton cycles
in~\cite{W6}. Even though both proofs in~\cite{W7} and~\cite{W6} are
based on a similar counting trick, the analysis for the undirected
version is much more complicated. The proof for the directed
Hamilton cycles was redone by Frieze and Suen~\cite{FS}, probably
unaware of the existing work of Wrignt, using basically the same
approach. In~\cite{J3}, Janson reproved the same result for both the
undirected and directed versions, using the method of graph
decomposition and projection. In this paper, we present a much shorter proof, using our method, for both the directed
and undirected versions.

We also present two new results: the distribution of the number of triangle free subgraphs in Section~\ref{sec:Hfree}, and the distribution
of the number of triangle-factors (the spanning subgraphs isomorphic to a collection of vertex disjoint
triangles) in Section~\ref{sec:dsjTriangle}. Their distributions are determined by verifying the
conditions given in the theorems in Section~\ref{sec:general}. In Theorem~\ref{t:Gnm}, we
 state a general approach for proving
concentration of any large subgraphs in $\G(n,m)$. The proofs in~\cite[Theorems 2.3 and 2.4]{G4} implicitly follow
the approach as described in Theorem~\ref{t:Gnm}, though in~\cite{G4} the setting is only for examining the $d$-factors.
The idea of the proof of Theorem~\ref{t:Gnp} is essentially the same as the
proof of~\cite[Theorem 6]{J3} (with necessary modifications), except that we relax some constraints in~\cite[Theorem 6]{J3} so that it is applicable to the study of a larger class of subgraphs.
The proofs of both Theorems~\ref{t:Gnm} and~\ref{t:Gnp}
are presented in Section~\ref{sec:proofs}.  

\section{A general approach}
\lab{sec:general}

Let $\S$ denote a set of vertex-labelled graphs on a set $S=[n]$ of
$n$ vertices. 
For two graphs $H_1$ and $H_2$ both on vertex set $S$, let $H_1\cap
H_2$ ($H_1\cup H_2$) denote the set of edges contained in both
(either of) $H_1$ and $H_2$. For any integer $j\ge 0$, let
$F_{j}(\S)$ denote the set of ordered pairs $(H_1,H_2)\in
\S\times\S$ such that $| H_1 \cap H_2 |=j$. Let
$f_j=f_j(\S)=|F_j(\S)|$ and let $r_j=f_j/f_{j-1}$ for any $j\ge 1$,
as long as $f_{j-1}\neq 0$. Let $X_n=X_n(\S)$ denote the number of
members of $\S$ that are contained in a random graph ($\G(n,p)$ or
$\G(n,m)$, defined on the same vertex set $S$) as (spanning)
subgraphs.  Here $S$, $p$ and $m$ refer to sequences $(S(n))_{n\ge
1}$, $(p(n))_{n\ge 1}$ and $(m(n))_{n\ge 1}$. Assume every graph in
$\S$ has the same number $h=h(n)$ of edges. Let
$N(n)=\binom{n}{2}$. We drop $n$ from all these notations when there
is no confusion. All asymptotics in this paper refer to
$n\to\infty$. For any real $x$ and any integer $\ell\ge 0$, define
the $\ell$-th falling factorial $[x]_{\ell}$ to be
$\prod_{i=0}^{\ell-1}(x-i)$. Let
\begin{eqnarray}
\mu_n&=&|\S| \binom{N-h}{m-h}\Big/\binom{N}{m},\ \ \
\la_n=|\S|p^h.\lab{mu-la}
\end{eqnarray}
Clearly,
$$
\ex_{\G(n,m)} X_n=\mu_n,\ \ \ \ex_{\G(n,p)}X_n=\la_n.
$$
A simplification of $\mu_n$ (readers can also refer to
Lemma~\ref{l:cal} by taking $\ell=h$) gives
\begin{eqnarray}
\ex_{\G(n,m)}
X_n&=&|\S|\cdot\frac{\binom{N-h}{m-h}}{\binom{N}{m}}=|\S|\cdot\frac{[m]_{h}}{[N]_h}=|\S|(m/N)^h\exp\left(-\frac{N-m}{mN}\frac{h^2}{2}+O(h^3/m^2)\right).\lab{mu}
\end{eqnarray}
\begin{thm}\lab{t:Gnm} Let $\mu_n$ be defined as in~\eqn{mu-la}.
Assume that $h^3=o(m^2)$, $h^2=\Omega(m)$ and for $\rho(n)=h^2/m$ and some function
$\gamma(n)$, the following conditions hold:
\begin{description}

\item{(a)} for all $K>0$ and for all $1\le j\le K\rho(n)$,
$$
r_j=\frac{h^2}{Nj}(1+o(m/h^2));
$$
\item{(b)} $r_j\le m/2N$ for all $4\rho(n)\le j\le \gamma(n)$;
\item{(c)} $t(n):=\sum_{j>\gamma(n)}f_j=o(\mu_n|\S|)$.
\end{description}
 Then in $\G(n,m)$,
$$
X_n/\ex_{\G(n,m)}(X_n)\xrightarrow{p} 1,
$$
as $n\to\infty$.

\end{thm}

\no {\bf Remark}: The ratio $r_j$ in condition (a) looks quite
restrictive. However, as we will see in the next section, this ratio
appears naturally if the edges in $\S$ are distributed randomly (see
examples in Sections~\ref{sec:trivial} and~\ref{sec:trivial2}). In
some cases, for instance, if we take $\S$ to be the set of graphs
isomorphic to a given unlabelled graph on $n$ vertices, the edges in
$\S$ are likely to still distribute in some kind of ``random-like''
way and thus having $r_j$ as expressed in condition (a) is expected.
If we are lucky, we might have condition (b) satisfied for
$\gamma(n)=h$. Then $t_n=0$ and condition (c) is satisfied trivially. See the example in Section~\ref{sec:Hfree}. But
usually this is not the case, as the sequence $r_j$ might decrease
first and increase at its tail. Normally, in these cases, condition
(c) is not difficult to verify. See examples in
Sections~\ref{sec:hamilton}
and~\ref{sec:dsjTriangle}. \ss

Theorem~\ref{t:Gnm} and its proof also gives the following
proposition.

\begin{prop} Assume all conditions (a)--(c) of Theorem~\ref{t:Gnm} are
satisfied with $m=N$. Then,
 for all $j=O(h^2/m)$,
 $$
 f_j(n)\sim |\S|^2\exp(-h^2/N)(h^2/N)^j/j!.
 $$

\end{prop}

The following theorem gives conditions under which $X_n$ will be asymptotically log-normally distributed in $\G(n,p)$ if all conditions in Theorem~\ref{t:Gnm} are satisfied by taking $m=pN$.

\begin{thm}\lab{t:Gnp}
Assume $h^3=o(p^2n^4)$. Let $\beta_n=h\sqrt{(1-p)/pN}$ and $\la_n$
as defined in~\eqn{mu-la}. Assume further that
$\liminf_{n\to\infty}\beta_n>0$. If for all $m=pN+O(\sqrt{pN})$,
$X_n/\ex_{\G(n,m)}(X_n)\xrightarrow{p} 1$, then
 $$
\frac{\ln(e^{\beta_n^2/2} X_n/\la_n)}{\beta_n}\xrightarrow{d}
\N(0,1), \ \ \mbox{as}\ n\to\infty,
$$
where $\N(0,1)$ is the standard normal distribution.
\end{thm}

Combining Theorems~\ref{t:Gnm} and~\ref{t:Gnp}, we immediately have the following corollary.
\begin{cor}\lab{c:combine}
Assume $h^3=o(p^2n^4)$ and $h^2=\Omega(pn^2)$. Let $\beta_n=h\sqrt{(1-p)/pN}$ and $\la_n$
as defined in~\eqn{mu-la}. Assume further that
$\liminf_{n\to\infty}\beta_n>0$. If for all $m=pN+O(\sqrt{pN})$, conditions (a)--(c) in Theorem~\ref{t:Gnm} are satisfied, then
$$
\frac{\ln(e^{\beta_n^2/2} X_n/\la_n)}{\beta_n}\xrightarrow{d}
\N(0,1), \ \ \mbox{as}\ n\to\infty,
$$
where $\N(0,1)$ is the standard normal distribution.
\end{cor}

Hence, in order to prove that a subgraph count has a log-normal distribution in $\G(n,p)$, it is enough to
check conditions (a)--(c) in Theorem~\ref{t:Gnm} by taking $m=pN+O(\sqrt{pN})$ if the value of $p$ and the number of edges in the subgraph $h$ satisfy the hypotheses in Corollary~\ref{c:combine}.
This method is particularly powerful if we can estimate $r_j$ without
knowing $f_j$. This is usually the case when we apply the switching method developed by McKay~\cite{M2}.
As we will see in the later examples, our method is easy to be applied by making
extensive use of the switching method.

We can generalise the results to random digraphs. Define $\D(n,m)$
to be the random digraph on $n$ vertices with $m$ directed edges
chosen uniformly at random from the $2N$ ordered pairs of vertices.
Define $\D(n,p)$ to be the random digraph on $n$ vertices, which
includes every directed edge independently with probability $p$. In
this paper, we again define $\D(n,m)$ and $\D(n,p)$ on the vertex
set $S$. With almost the same proofs of Theorems~\ref{t:Gnm}
and~\ref{t:Gnp} we have the following theorem.

\begin{thm}\lab{t:digraph} The same conclusions of Theorems~\ref{t:Gnm} and~\ref{t:Gnp}
hold if we replace $\G(n,m)$, $\G(n,p)$, $N$ by $\D(n,m)$, $\D(n,p)$
and $2N$.

\end{thm}

\section{Two trivial examples}

The purpose of this section is to provide simple demonstrations of our method and to convince readers that the behavior of the ratio $r_j$ as given
in Theorem~\ref{t:Gnm} (a) shall be well expected.

\lab{sec:examples}

\subsection{The first trivial example}\lab{sec:trivial} Take $\S_1$ to be the set of all graphs on vertex
set $S$ with $h$ edges. Then $|\S_1|=\binom{N}{h}$. The conclusion
of Theorem~\ref{t:Gnm} should hold trivially in this case as
$X_n(\S_1)$ is constant in $\G(n,m)$ (depending only on $m$ and $h$).
Nevertheless we verify conditions (a) and (b), also for later use in
the next section. For all $0\le j\le h$,
$$
f_j=\binom{N}{j}\binom{N-j}{h-j}\binom{N-h}{h-j}.
$$
Then for all $1\le j\le h$,
$$
r_j=\frac{(h-j+1)^2}{j(N-2h+j)}=\frac{h^2}{jN}(1+O(j/h+h/n^2)).
$$
This verifies conditions (a) and (b) (for $\gamma(n)=h$).

\subsection{Another trivial example}\lab{sec:trivial2}

Let $0<\hat p<1$. Consider the set of graphs $\S_2$ that is obtained
by including each element in $\S_1$ independently with probability
$\hat p$. Note here that $\S_2$ itself is a random variable. Then we have the following.

\begin{thm}\lab{t:ranSet}
Assume $0<\hat p\le 1$, $0<p<1$ are reals and $h$ is an
integer that satisfy $m=p\binom{n}{2}$, $h^3=o(m^2)$, $h^2=\Omega(m)$, $m^2\hat p^2
N^h>> h^{3h+4}\ln n$. 
Let $\mu_n(\S_2)$ and $\lambda_n(\S_2)$ be defined as in~\eqn{mu-la} and let
$\beta_n=h\sqrt{(1-p)/pN}$. Then 
$X_n(\S_2)/\mu_n(\S_2)\xrightarrow{p} 1$  in $G(n,m)$, and
$$
\frac{\ln(e^{\beta_n^2/2}X_n(\S_2)/\lambda_n(\S_2))}{\beta_n}\xrightarrow{d}
\N(0,1),\ \ \mbox{in}\ \G(n,p),
$$ provided $\liminf_{n\to\infty}\beta_n>0$.
\end{thm}

\proof By the definition of $\S_2$, we have $s_2=|\S_2|\sim
\B(\binom{N}{h},\hat p)$ and for all $0\le j<h$, $f_j\sim \B(M_j,\hat p^2)$, $f_h\sim\B(M_h,\hat p)$, where
$M_j=\binom{N}{j}\binom{N-j}{h-j}\binom{N-h}{h-j}$.
Let $A_n$ denote the family of $\S_2$ which satisfies
$$
\forall\ 0\le j\le h-1,\ f_j=\left(1+o(m/h^2)\right)\hat p^2 M_j,\quad f_h<2M_h\hat p,\quad s_2>\binom{N}{h}\hat p/2. 
$$

 The Chernoff
bound gives that
$$
\pr(|f_j-\hat p^2 M_j|>2\sqrt{3(\ln n)\hat p^2M_j})<\exp(- 3\ln
n)=n^{-3},\ \ \forall\ 0\le j\le h-1,
$$
and
$$
\pr(f_h<2M_h\hat p)=1-o(1),\quad \pr\left(s_2>\binom{N}{h}\hat p/2\right)=1-o(1).
$$
Therefore, with probability at least $1-hn^{-3}-o(1)=1-o(1)$, for
all $0\le j\le h-1$,
$$
f_j=\left(1+O\left(\sqrt{\ln n/\hat p^2M_j}\right)\right)\hat p^2 M_j.
$$
Note that for all $j$, $M_j>[N]_h/(h!)^3>(N/h^3)^h$ and $\hat p$
satisfies
$$
\hat p^2>>\frac{h^{4+3h}\ln n}{m^2N^h}.
$$
Thus, a.a.s.\ for all $0\le j\le h-1$,
$$
f_j=\left(1+o(m/h^2)\right)\hat p^2 M_j,
$$
i.e.\ $\pr(S_2\in A_n)=1-o(1)$. For every $S_2\in A_n$,
by the calculations in Section~\ref{sec:trivial}, both
conditions (a) and (b) (for $\gamma(n)=h-1$) are satisfied whereas condition (c) can be easily verified by noting that
$
t_n=f_h< 2M_h\hat p=o(\mu_n|\S_2|)$. The theorem
thereby follows.\qed

\medskip

The following is a corollary of Theorem~\ref{t:ranSet} by letting
$\hat p=1/2$. Here $\S'_2$ are no longer random variables. We may consider $\S'_2$ as elements in $A_n$ in the proof of Theorem~\ref{t:ranSet}.

\begin{cor}\lab{c:ranSet}
Assume $0<p<1$ is a real and $h$ is an integer that
satisfy $m=p\binom{n}{2}$, $h^3=o(m^2)$, $h^2=\Omega(m)$, $m^2 N^h>> h^{3h+4}\ln n$.
Then for almost all subsets $\S'_2$ of $\S_1$, the same conclusions
of Theorem~\ref{t:ranSet} hold when $\S_2$ is
replaced by $\S'_2$.
\end{cor}

\section{A new approach -- Hamilton cycles}
\lab{sec:hamilton} The most interesting examples of $\S$ are perhaps
taking $\S$ as the set of graphs that are isomorphic to a given
unlabelled graph $H$ on a set of $n$ vertices. In this section, we investigate the number of Hamilton cycles. In literature, computing
the second moment of the number of Hamilton cycles involves heavy analysis, as done by Wright~\cite{W6,W7}, using the inclusion and exclusion and
some recursive functions, and by Janson~\cite{J3}, using the graph decompostion and projection. Here, we present a new and much shorter
proof.

Let $H$ ($H'$) be a cycle (directed
cycle) with length $n$ and $\S_3$ ($\S'_3$) to be the set of graphs
(directed graphs) on $S$ that are isomorphic to $H$ ($H'$). Thus,
$X_n(\S_3)$ and $X_n(\S'_3)$ count the numbers of undirected and
directed Hamilton cycles respectively. It is well known that
\begin{equation}
|\S_3|=(n-1)!/2,\ \mbox{and} \ |\S'_3|=(n-1)!.\lab{s5}
\end{equation}

We have the following theorem for the undirected version.

\begin{thm}\lab{t:hamilton}Let $0<p<1$ be a real and $0<m<N$ an integer satisfying $m=pN$ and $p>>n^{-1/2}$. Let $X_n$ denote the number of Hamilton cycles in $\G(n,m)$ (or $\G(n,p)$).
 Let $\mu_n=\ex_{\G(n,m)}X_n$ and
let $\lambda_n=\ex_{\G(n,p)}X_n$. Then $ X_n/\mu_n\xrightarrow{p} 1$
in $G(n,m)$.  Assume further that $\limsup_{n\to\infty} p(n)<1$,
then
$$
\frac{\ln(e^{\beta_n^2/2}X_n/\lambda_n)}{\beta_n}\xrightarrow{d}
\N(0,1),\ \ \mbox{in}\ \G(n,p),
$$
where $\beta_n=\sqrt{2(1-p)/p}$.
\end{thm}

\proof We define two
switching operations as follows.

 \no {\em $h$-switching}: Choose an edge $xy\in G_1\cap
G_2$. Then choose edges $x_1y_1\in G_1\setminus G_2$, $x_2y_2\in
G_2\setminus G_1$ such that $xyx_1y_1$ and $xyx_2y_2$ are in a
cyclic order in $G_1$ and $G_2$ respectively. Replace $xy$ and
$x_1y_1$ by $xx_1$ and $yy_1$ in $G_1$, and replace $xy$ and
$x_2y_2$ by $xx_2$ and $yy_2$ in $G_2$. The $h$-switching is
applicable if and only if
\begin{description}
\item{(a)} the six vertices $x$, $y$, $x_i$ and $y_i$ for $i=1,2$
are all distinct;
\item{(b)} the edges $xx_1$ and $yy_1$ are not in $G_2$ and the
edges $xx_2$ and $yy_2$ are not in $G_1$.
\end{description}

\no{\em inverse $h$-switching}: Choose a pair of vertices $\{x,y\}$
such that $xy\notin G_1\cup G_2$. For $i=1,2$, choose $x_i$ and
$y_i$ such that $xx_i\in G_i$ and $yy_i\in G_i$ and $xx_iyy_i$ is in
a cyclic order in $G_i$. The inverse $h$-switching replaces $xx_i$
and $yy_i$ by $xy$ and $x_iy_i$ in $G_i$ for $i=1,2$. The operation
is applicable if and only if
\begin{description}
\item{(a')} the six vertices $x$, $y$, $x_i$ and $y_i$ for $i=1,2$
are all distinct;
\item{(b')} the edges $xx_i$ and $yy_i$ are not in $G_1\cap G_2$ for
$i=1,2$;
\item{(c')} $x_1y_1\notin G_2$ and $x_2y_2\notin G_1$.
\end{description}

For $g\in F_j$, let $N(g)$ be the number of $h$-switchings
applicable on $g$. There are $2j$ ways to choose and label the end
vertices of the edge $xy\in G_1\cap G_2$. For any chosen $xy$, there
are $n-j+O(1)$ ways to choose and label the end vertices of the edge
$x_iy_i\in G_i$, where $j+O(1)$ accounts for the case that $x_iy_i
\in G_1\cap G_2$ and the case that condition (a) is violated. Thus,
a rough estimation of $N(g)$ is $2j(n-j+O(1))^2$. The only miscounts
are those $xy$ and $x_iy_i$ such that condition (b) is violated.
Clearly, the miscount due to the violation of condition (b) is
$O(jn)$ because for any chosen $xy$, there are exactly two choices
for $x_1y_1$ (equivalently $x_2y_2$), such that either $xx_1$ or
$yy_1$ is in $G_2$ (equivalently, either $xx_2$ or $yy_2$ is in
$G_1$). Thus, $ N(g)=2jn^2(1-j/n+O(n^{-1}))^2$.

On the other hand, for $g'\in F_{j-1}$, let $N'(g')$ denote the
number of inverse $h$-switchings applicable on $g'$. There are
$n^2-O(n)$ ways to choose and label vertices $x$ and $y$ such that
$xy\notin G_1\cup G_2$. For any chosen $xy$, there are two ways to
choose $x_i$ and $y_i$ from $G_i$ for $i=1,2$ respectively, such
that $xx_i$, $yy_i\in G_i$ and $xx_iyy_i$ is in a cyclic order in
$G_i$. Thus, $N'(g')$ is approximately $4(n^2-O(n))$. The only
miscounts are those choices that violate conditions (a') or (b') or
(c'). There are only $O(n)$ choices of $xy$ so that (a') or (c') can
possibly be violated, and there are only $O(jn)$ choices of $xy$ so
that (b') can possibly be violated. Therefore,
$N'(g')=4n^2(1+O(j/n))$.

Hence for all $1\le j\le n/2$,
$$
r_j=\frac{4n^2}{2jn^2}(1+O(j/n))=\frac{2}{j}(1+O(j/n)),
$$
from which we can easily verify Theorem~\ref{t:Gnm} (a), (b) (for
$\gamma(n)=n/2$). 
 The proof will
be completed by verifying condition (c). Let $G$ be a Hamilton
cycle, and let $\kappa_j(G)$ denote the number of Hamilton cycles
that share at least $j$ edges with $G$. There are $\binom{n}{j}$
ways to choose $j$ edges from $G$. These chosen edges form $r\le j$
disjoint paths. Contract each path into a special vertex. The total
number of vertices including these special vertices is then $n-j$.
There are $(n-j-1)!/2$ Hamilton cycles on these vertices. For every
such Hamilton cycles, expand each special vertex by its
corresponding path (there are two ways to expand each special
vertex). Then each expanded Hamilton cycle corresponds to a Hamilton
cycle that shares at least $j$ edges with $G$. Thus, for every $G$,
$$
\kappa_j(G)\le \binom{n}{j}\frac{(n-j-1)!}{2}\cdot 2^j<n! 2^j/j!.
$$
It is then straightforward to verify that
$$
\sum_{j\ge n/2}f_j\le |\S_3| n! 2^{n/2}/(n/2)!=o(|\S_3|\mu_n).\qed
$$

The same proof, with only slight modification
of the switchings that cope with directed edges, works for the
directed version (Theorem~\ref{t:dirHamilton}).

\begin{thm} \lab{t:dirHamilton} If all assumptions with $N$, $\G(n,p)$ and $\G(n,m)$ replaced
by $2N$, $\D(n,p)$ and $\D(n,m)$ in Theorem~\ref{t:hamilton}  hold,
then the same conclusion of Theorem~\ref{t:hamilton} holds (for
$\beta_n=\sqrt{(1-p)/p}$ by the definition of $\beta_n$ in
Theorem~\ref{t:Gnp}).
\end{thm}

\section{triangle-free subgraphs}
\lab{sec:Hfree}   In this section, we consider another example where
$\S_4$ is the set of all triangle-free graphs on $S$ with $h$ edges
and maximum degree at most $\Delta=\Theta(h^{1/3})$. Then $X_n(\S_4)$ counts the
number of triangle-free subgraphs with $h$ edges and maximum degree at most $\Delta$.

\begin{thm}\lab{t:triangle} Let $0<p<1$ be a real and $0<m<N$ an integer satisfying $m=pN$, $h^2=\Omega(m)$, $h^3=o(m^2)$ (or equivalently
$h^3=o(p^2n^4)$) and $h^{8/3}=o(pn^3)$. Assume $\Delta=\Theta(h^{1/3})$ is an integer. Let $X_n$ denote the number of
triangle-free subgraphs with $h$ edges and maximum degree
at most $\Delta$. Let $\mu_n$ and
 $\lambda_n$ be defined as in~\eqn{mu-la} and let $\beta_n=h\sqrt{(1-p)/pN}$. Then
$ X_n/\mu_n\xrightarrow{p} 1$ in $G(n,m)$, and
$$
\frac{\ln(e^{\beta_n^2/2}X_n/\lambda_n)}{\beta_n}\xrightarrow{d}
\N(0,1),\ \ \mbox{in}\ \G(n,p),
$$provided
$\liminf_{n\to\infty}\beta_n>0$.
\end{thm}

\proof Recall that $F_j(\S_4)=\{(G_1,G_2)\in \S_4\times\S_4:
|G_1\cap G_2|=j\}$. Consider $j\ge 1$ and the classes $F_j(\S_4)$
and $F_{j-1}(\S_4)$. Let $K_n$ denote the complete graph on $S$. We
define two other switchings operating on $\S_4\times\S_4$ as
follows.

{\em $s_2$-switching}: Let $x$ be an edge in $G_1\cap G_2$. Choose
$y$ and $z$ from $K_n\setminus  G_1\cup G_2 $, such that
\begin{description}
\item{(a)} $G_1\cup y$ and $G_2\cup z$ are triangle-free;
\item{(b)} $y$ ($z$) is not incident with a vertex with degree equal to
$\Delta$ in $G_1$ ($G_2$).
\end{description}
Replace $x$ by $y$ in $G_1$ and replace $x$ by $z$ in $G_2$.

{\em inverse $s_2$-switching}: Let $x$ be an edge in $K_n\setminus
 G_1\cup G_2$ such that
\begin{description}
\item{(a')}  $G_1\cup x$ and $G_2\cup x$ are triangle-free;
\item{(b')} In both $G_1$ and $G_2$, $x$ is not incident with a vertex with
degree equal to $\Delta$.
\end{description}
Let $y\in  G_1 \setminus  G_2 $ and $z\in
 G_2 \setminus  G_1 $. Replace $y$ by $x$ in $G_1$ and replace $z$
by $x$ in $G_2$.

Clearly, an $s_2$-switching converts an element $g\in F_j(\S_4)$ to
an element $g'\in F_{j-1}(\S_4)$ and an inverse $s_2$-switching
converts an element $g'\in F_{j-1}(\S_4)$ to an element $g\in
F_j(\S_4)$ for some $j\ge 1$. For any $g\in F_j(\S_4)$, let $N(g)$
denote the number of $s$-switchings that are applicable on $g$. Note
that in both $G_1$ and $G_2$, the number of vertices with degree
equal to $\Delta$ is $O(h^{2/3})$. There are $j$ ways to
choose $x$. Given $x$, the number of ways to choose $y$ and $z$ is
$N-O(h+nh^{2/3}+T_1(g))$ and $N-O(h+nh^{2/3}+T_2(g))$
respectively, where $T_i(g)$ denotes the number of $2$-paths in
$G_i$, and $O(nh^{2/3})$ bounds the number of forbidden choices such that $y$ (or $z$) is incident to a vertex with degree equal to $\Delta$.  Let $T(g)=\max\{T_1(g),T_2(g)\}$. Clearly
$T(g)=O(n\Delta^2)=O(nh^{2/3})$. So $N(g)=j(N-O(nh^{2/3}))^2$.
Then
$N(g)=jN^2(1+O(h^{2/3}/n))$. For any $g'\in F_{j-1}(\S_4)$, let
$N'(g')$ denote the number of inverse $s_2'$-switchings applicable
on $g'$. Then
$N'(g')=(N-O((2h-j+1)+nh^{2/3}+T(g')))(h-j+1)^2=Nh^2(1+O(h^{2/3}/n+j/h))$.
Since $\sum_{g\in F_j(\S_4)} N(g)=\sum_{g'\in F_{j-1}(\S_4)}N'(g')$,
we have that for all $j\ge 1$,
\begin{equation}
r_j=\frac{Nh^2}{jN^2}(1+O(h^{2/3}/n+j/h))=\frac{h^2}{jN}(1+o(m/h^2)+O(j/h)).\lab{ratio5}
\end{equation}
Note that $O(h^{2/3}/n)=o(m/h^2)$ because  $h^{8/3}=o(pn^3)$. Next we
verify conditions (a) and (b) of Theorem~\ref{t:Gnm}. For all
$j=O(h^2/m)$, $j/h=O(h/m)=o(m/h^2)$ since $h^3=o(m^2)$. Thus
$$
r_j=\frac{h^2}{jN}(1+o(m/h^2)),
$$
which verifies condition (a). By~\eqn{ratio5}, for all $j\ge
3h^2/m$,
$$
r_j=\frac{h^2}{jN}(1+o(1))+O(h/N)\le \frac{m}{2N},
$$
which verifies condition (b) (for $\gamma(n)=h$). \qed

\section{Triangle-factors}
\lab{sec:dsjTriangle}
Given a graph $G$ on $n$ vertices where $n$ is a multiple of $3$, a subgraph of $G$ consisting of $n/3$ vertex disjoint triangles is called
a triangle-factor of $G$.
 In this section, we assume $n\equiv 0\
(\mbox{mod}\ 3)$ and consider $H$ ($H'$) to be the unlabelled graph on $n$ vertices consisting of $n/3$ vertex disjoint triangles (directed triangles). Let $\S_5$ ($\S'_5$)
denote the set of graphs on $S$ that are isomorphic to $H$ ($H'$). Then $X_n(\S_5)$ counts the number of triangle-factors and
\begin{equation}
|\S_5|=\frac{n!}{6^{n/3}(n/3)!},\ \ \ |\S'_5|=\frac{n!}{3^{n/3}(n/3)!}.\lab{s6}
\end{equation}

The following theorem determines the limiting distribution of
$X_n=X_n(\S_6)$.

\begin{thm}\lab{t:dsjTriangle}
Let $0<p<1$ be a real and $0<m<N$ an integer satisfying $m=pN$ and
$\liminf_{n\to\infty}p(n)>0$. Let $X_n$ denote the number of
subgraphs that are isomorphic to a set of $n/3$ vertex disjoint
triangles.
 Let $\mu_n=\ex_{\G(n,m)}X_n$ and
let $\lambda_n=\ex_{\G(n,p)}X_n$. Then $X_n/\mu_n\xrightarrow{p} 1$
in $\G(n,m)$.  Assume further that $\limsup_{n\to\infty} p(n)<1$,
then
$$
\frac{\ln(e^{\beta_n^2/2}X_n/\lambda_n)}{\beta_n}\xrightarrow{d}
\N(0,1),\ \ \mbox{in}\ \G(n,p),
$$
where $\beta_n=\sqrt{2(1-p)/p}$.
\end{thm}

\no {\bf Remark}: Indeed, the condition of $\liminf_{n\to\infty}
p(n)>0$ can be replaced by $p(n)\ge n^{-\delta}$, for some small
constant $\delta$. For instance, we checked that $\delta=1/16$ works
and there is still room for further improvement. However,
$p>>n^{-1/2}$ does not seem to be sufficient. For the purpose of a
cleaner presentation, we only consider $\liminf_{n\to\infty} p(n)>0$
in the proof. For readers who are interested in improving the
condition of $p$, we give quite tight bounds in
Lemmas~\ref{l:single} and~\ref{l:triple}, and we also point out here
that there is plenty of room in the proofs of Lemma~\ref{l:largeJ}
and Theorem~\ref{t:dsjTriangle} to improve the range of $p$.\ss

Almost the same proof of the previous theorem, with slight modifications of the switchings defined in the proof of
Theorem~\ref{t:dsjTriangle}, concerning the directions of edges, yields the following corresponding theorem for
the number of directed triangle-factors.

\begin{thm}\lab{t:dirTriangle} If all assumptions with $N$, $\G(n,p)$ and $\G(n,m)$
replaced by $2N$, $\D(n,p)$ and $\D(n,m)$ in
Theorem~\ref{t:dsjTriangle}hold, then the same conclusion of
Theorem~\ref{t:dsjTriangle} holds (for $\beta_n=\sqrt{(1-p)/p}$ by
the definition of $\beta_n$ in Theorem~\ref{t:Gnp}).
\end{thm}

 For
any $(G_1,G_2)\in \S_6\times\S_6$, the edges in $G_1$ and $G_2$ can
intersect in two ways. We say $e\in G_1\cap G_2$ is of type 1 if the
triangles $T_i\in G_i$ with $e\in T_i$ for $i=1,2$ are distinct. We
say $e$ is of type 2 if $T_1$ and $T_2$ are on the same vertex set.

Let $F_{\ell,t}$ denote the set of $(G_1,G_2)\in \S_6\times\S_6$
such that number of edges in $G_1\cap G_2$ of type 1 and 2 is $\ell$
and $t$ respectively. Clearly $F_{\ell,t}$ is non-empty only if $t$
is a multiple of $3$. Clearly $F_j(\S_6)=\cup_{k}F_{j-3k,3k}$. Let
$f_{\ell,t}=|F_{\ell,t}|$. Then $f_j=\sum_{k=0}^{\lfloor
j/3\rfloor}f_{j-3k,3k}$.

\begin{lemma}\lab{l:single} For any $t\ge 0$ and $\ell\ge 1$ such that $n-4\ell-3t-1>0$ and $n-3\ell-3t-12>0$,
$$
\frac{2}{\ell}\frac{(n-4\ell-3t-1)^2}{(n-3\ell-3t)^2}\le\frac{f_{\ell,3t}}{f_{\ell-1,3t}}\le\frac{2}{\ell}\frac{(n-4\ell-3t+4)^2}{(n-3\ell-3t-12)^2}.
$$
\end{lemma}

\proof We define two switchings operating on $\S_6\times\S_6$ as
shown in Figure~\ref{f:t1}.

\no{\em $t_1$-switching}: Take an edge of type 1 in $G_1\cap G_2$
and label the end vertices $x$ and $y$. Let $u$ ($v$) be the vertex
that is adjacent to both $x$ and $y$ in $G_1$ ($G_2$). Take a
triangle $T_1$ ($T_2$) in $G_1$ ($G_2$) that is distinct from $xyu$
($xyv$) which does not contain any edge in $G_1\cap G_2$. Label the
vertices of $T_1$ ($T_2$) as $u_1u_2u_3$ ($v_1v_2v_3$). Replace
these four triangles in $G_1\cup G_2$ by $xuu_1$, $yu_2u_3\in G_1$
and $xvv_1$, $yv_2v_3\in G_2$. The $t_1$-switching is applicable
only if $v\notin T_1$, $u\notin T_2$ and $T_1\cap T_2=\emptyset$.
See Figure~\ref{f:t1}.\ss

\no{\em inverse $t_1$-switching}: A vertex $x$ is pure if both
triangles containing $x$ in $G_1$ and $G_2$ do not contain any edge
in $G_1\cap G_2$. Choose a pure vertex $x$ and label its neighbours
in $G_1$ ($G_2$) as $u$ and $u_1$ ($v$ and $v_1$). Then choose
another pure vertex $y$ that is distinct from $x$, $u_i$ and $v_i$
for $i=1,2$. Label the neighbours of $y$ in $G_1$ ($G_2$) as $u_2$
and $u_3$ ($v_2$ and $v_3$). Replace these four triangles under
consideration by $xyu$, $u_1u_2u_3\in G_1$ and $xyv$, $v_1v_2v_3\in
G_2$.\ss

  \begin{figure}[htb]
\vbox{\vskip .6cm
 \hbox{\centerline{\includegraphics[width=15cm]{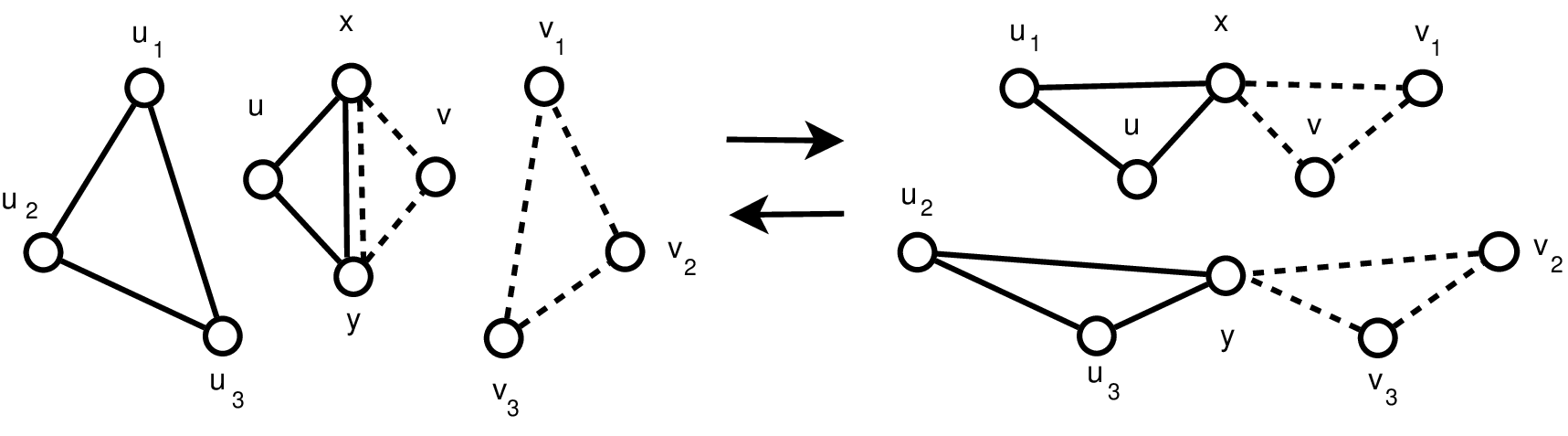}}}
\vskip .6cm \smallskip} \caption{\it  $t_1$-switching and its
inverse}

\lab{f:t1}

\end{figure}

For any $g=(G_1,G_2)\in F_{\ell,3t}$, let $N(g)$ be the number of
$t_1$-switchings that are applicable on $g$. Clearly $N(g)\le 2\ell
(6(n/3-(\ell+t)))^2$, as there are $2$ ways to label $x$ and $y$ for
a chosen edge from $G_1\cap G_2$, and in $G_1$ ($G_2$) there are at
most $n/3-(\ell+t)$ choices for the triangle $u_1u_2u_3$
($v_1v_2v_3$) and for each choice there are $6$ ways to label the
vertices. We also have
$$
N(g)\ge 2\ell\cdot 6(n/3-(\ell+t)-1)\cdot 6(n/3-(\ell+t)-4),
$$
because for any chosen $xy$, the number of triangles in $G_1$ which
contain no edges in $G_2$ and do not contain $v$ is at least
$n/3-(\ell+t)-1$, whereas given the triangle $u_1u_2u_3$, the number
of triangles in $G_1$ which contain no edges in $G_1$ and do not
contain any of $u$, $u_i$, $i=1,2,3$ is at least $n/3-(\ell+t)-4$.
  On the other hand, for any $g'=(G_1,G_2)\in
F_{\ell-1,3t}$, let $N'(g')$ be the number of inverse
$t_1$-switchings applicable on $g'$. The number of pure vertices is
exactly $n-4(\ell-1)-3t$. Hence the number of ways to choose $x$ is
$n-4(\ell-1)-3t$ and for any chosen $x$, the number of ways to label
$u$, $u_1$, $v$, $v_1$ is $4$. The number of ways to choose $y$ is
$n-4(\ell-1)-3t-\delta$, where $\delta$ counts the number of pure
vertices among $x$, $u$, $u_1$, $v$ and $v_1$. Therefore, $1\le
\delta\le 5$ always. Hence,
$$
\frac{16
 (n-4(\ell-1)-3t-5)^2}{2\ell\cdot(6(n/3-(\ell+t)))^2}\le\frac{f_{\ell,3t}}{f_{\ell-1,3t}}\le\frac{16(n-4(\ell-1)-3t)^2}{2\ell\cdot
36(n/3-(\ell+t)-4)^2}.\qed
$$

\begin{lemma}\lab{l:triple} For any $\ell\ge 0$ and $t\ge 1$,
$$
\frac{f_{\ell,3t}}{f_{\ell,3(t-1)}}=\frac{32(n-4\ell-3t)^3}{3(n-3\ell-3t)^4}(1+O(1/(n-4\ell-3t))).
$$
\end{lemma}

\proof We define another two switching operations on
$\S_6\times\S_6$ as shown in Figure~\ref{f:t2}.

\no {\em $t_2$-switching}: Let $xyz$ be a triangle that is contained
in both $G_1$ and $G_2$. Take two distinct triangles from $G_1$
($G_2$) which do not contain any edge in $G_1\cap G_2$ and label the
end vertices as $x_1y_1z_1$ and $x_2y_2z_2$ ($x'_1y'_1z'_1$ and
$x'_2y'_2z'_2$) respectively. Replace the six triangles under
consideration by $aa_1a_2\in G_1$ and $aa'_1a'_2\in G_2$, where
$a\in\{x,y,z\}$. This switching is applicable only if all these
fifteen vertices $a$, $a_i$, $a'_i$ for $a\in\{x,y,z\}$ and $i=1,2$
are distinct.\ss

\no{\em inverse $t_2$-switching}: Recall from the definition of
inverse $t_1$-switching that a vertex $x$ is pure if both triangles
containing $x$ in $G_1$ and $G_2$ do not contain any edge in
$G_1\cap G_2$. Choose three pure vertices $a$, $a\in\{x,y,z\}$ and
label the neighbours of $a$ in $G_1$ ($G_2$) by $a_1$ and $a_2$
($a'_1$ and $a'_2$). The inverse $t_2$-switching replaces the six
triangles under consideration by $xyz$, $x_iy_iz_i\in G_1$ for
$i=1,2$ and $xyz$, $x'_iy'_iz'_i\in G_2$ for $i=1,2$. This switching
is applicable only if all these fifteen vertices $a$, $a_i$, $a'_i$
for $a\in\{x,y,z\}$ and $i=1,2$ are distinct.

  \begin{figure}[htb]
\vbox{\vskip .6cm
 \hbox{\centerline{\includegraphics[width=15cm]{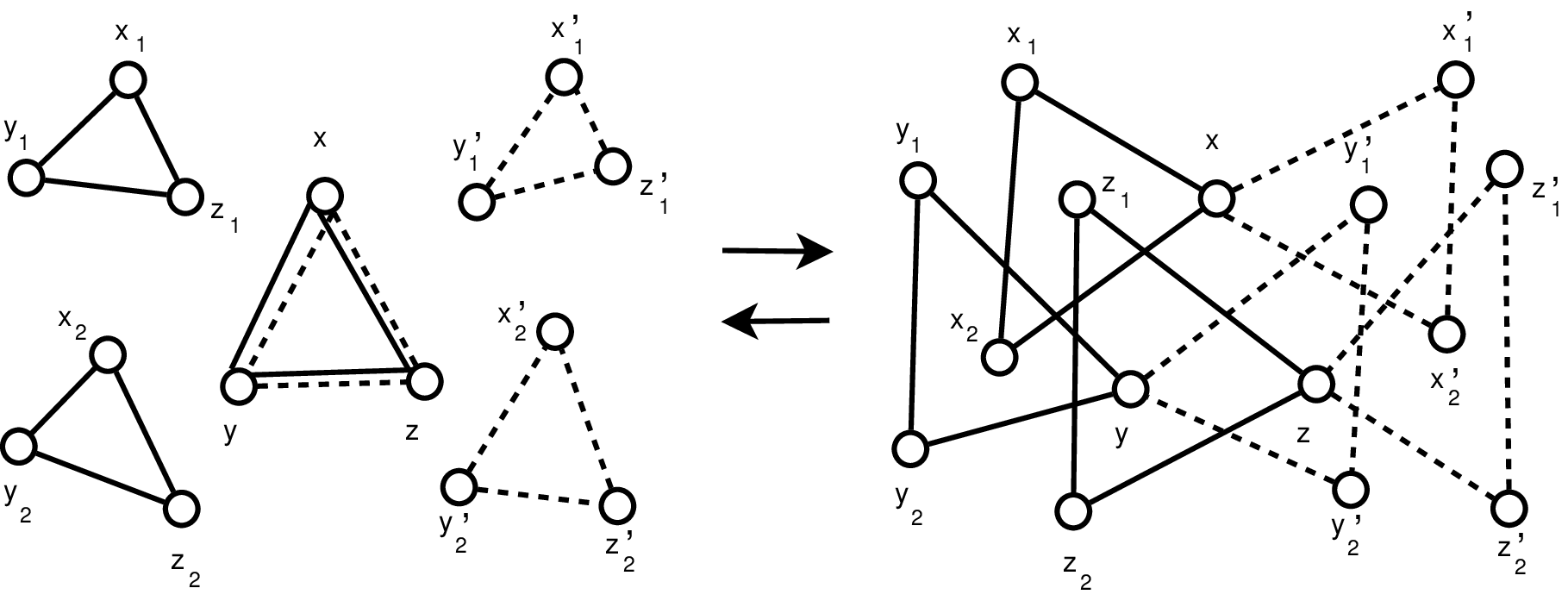}}}
\vskip .6cm \smallskip} \caption{\it  $t_2$-switching and its
inverse}

\lab{f:t2}

\end{figure}

For any $g\in F_{\ell,3t}$ and $g'\in F_{\ell,3t-3}$, define $N(g)$
and $N'(g')$ the same way as in the proof of Lemma~\ref{l:single}.
Following an analogous argument of Lemma~\ref{l:single}, it is not
hard to show that
\begin{eqnarray*}
&&6t\cdot6^2 \binom{n/3-(\ell+t)-6}{2} ^2\le N(g)\le 6t\cdot6^2 \binom{n/3-(\ell+t)}{2} ^2 \\
&&(4(n-4(\ell-1)-3t-10))^3\le N'(g')\le (4(n-4(\ell-1)-3t))^3.
\end{eqnarray*}
  Thus,
$$\frac{32(n-4\ell-3t-6)^3}{3(n-3\ell-3t)^4}\le
\frac{f_{\ell,3t}}{f_{\ell,3(t-1)}}\le
\frac{32(n-4\ell-3t+4)^3}{3(n-3\ell-3t-21)^4}.\qed
$$

\begin{cor}\lab{c:triple} For all $j=o(n)$,
$$
\frac{f_{j-3k-3,3k+3}}{f_{j-3k,3k}}\sim\frac{4[j-3k-1]_3}{3n}.
$$
\end{cor}

\proof This follows by Lemmas~\ref{l:single} and~\ref{l:triple} and
$$
\frac{f_{j-3k-3,3k+3}}{f_{j-3k,3k}}=\frac{f_{j-3k-3,3k+3}}{f_{j-3k-3,3k}}\prod_{i=0}^2\frac{f_{j-3k-i-1,3k}}{f_{j-3k-i,3k}}.\qed
$$

\begin{lemma}\lab{l:largeJ} Assume $\liminf_{n\to\infty}p(n)>0$. Let $\gamma(n)=n/\ln\ln n$. Then
$$
\sum_{j\ge \gamma(n)}f_j=o(|\S_6|\mu_n).
$$
\end{lemma}

\proof Let $G\in\S_6$ and let $\kappa_j(G)$ be the number of graphs
in $\S_6$ which shares at least $j$ edges with $G$. We estimate an
upper bound of $\kappa_j(G)$. Let $j=\ell+3t$ and we consider the
number of graphs $G'$ in $\S_6$ that shares at least $\ell$ and $3t$
edges of type 1 and 2 respectively with $G$. Then there are
$\binom{n/3}{t}$ ways to choose the $t$ triangles contained both in
$G$ and $G'$. Then there are $\binom{n/3-t}{\ell}3^{\ell}$ ways to
choose the $\ell$ triangles in $G$ that contain the $\ell$ edges of
type 1 and to locate these $\ell$ edges. Given these $\ell$ edges in
$G'$, there are at most $[n-3t-2\ell]_{\ell}$ ways to choose another
$\ell$ vertices to form the $\ell$ triangles in $G'$. Then there are
at most
$$
\frac{(n-3t-3\ell)!}{6^{n/3-t-\ell}(n/3-t-\ell)!}\le 9^n
n^{2(n/3-t-\ell)}
$$
ways to partition the remaining $n-3t-3\ell$ vertices into vertex
disjoint triangles in $G'$. Hence
$$
\kappa_j(G)\le \sum_{\ell}
\binom{n/3}{t}\binom{n/3-t}{\ell}3^{\ell}[n-3t-2\ell]_{\ell}9^n
n^{2(n/3-t-\ell)}\le n\cdot \max_{\ell}\{n^t
n^{2\ell}\ell^{-\ell}9^n n^{2(n/3-t-\ell)}\},
$$
where $t=(j-\ell)/3$. Thus,
$$
\ln(\kappa_j(G))\le \max_{\ell}\{ (2n/3-t)\ln
n-\ell\ln(\ell)\}+O(n).
$$
We consider only $j\ge\gamma(n)$. So the maximum is achieved at
$\ell=n^{1/3}$. Thus
$$
\ln(\kappa_j(G))\le  \frac{2n}{3}\ln n-\frac{j}{3}\ln n+O(n),
$$
We also have
$$
\ln\mu_n=n\ln p+\frac{2n}{3}\ln n+O(n).
$$
So
$$
\ln(\kappa_j(G))-\ln\mu_n\le-\frac{j}{3}\ln n-n\ln p+O(n)\to-\infty,
$$
as $n\to\infty$ since $\lim\inf_{n\to\infty} p(n)>0$, which
completes the proof of the lemma.\qed\ss

\remove{*****************************************
\begin{prop}\lab{p:largeT} Let $\alpha$ be a fixed integer. 
Then for any $j\ge n-\alpha$, $f_j=O(|\S_6|n^{3\alpha/2})$.
\end{prop}

\proof This assertion follows from two simple observations.

(a) We can write $j=\ell+t$, where $\ell$ and $t$ denote the number
of edges in $G_1\cap G_2$ of type 1 and 2 respectively. We observe
that $t+3\ell\le n$ since $t+3\ell$ counts the number of vertices
contained in triangles which contain edges in $G_1\cap G_2$. Thus,
$j\ge n-\alpha$ implies $t\ge n-3\alpha/2$.

(b) We have that $\sum_{\ell\le n-t}f_{\ell,t}\le
O(|\S_6|\binom{n}{t})$. This is because there are $\binom{n}{t}$
ways to choose the end vertices of the $t/3$ triangles contained in
both $G_1$ and $G_2$. The number of ways to partition these $t$
vertices into $t/3$ triangles is at most $|\S_6|$ and there are
finitely many $(A_1,A_2)$, where $A_i$ is a collection of $(n-t)/3$
vertex disjoint triangles among the remaining $n-t\le \alpha$
vertices.\qed\ss

**************************************************************************
**************************************************************************END
OF REMOVE}

\no {\bf Proof of Theorem~\ref{t:dsjTriangle}.\ } For any $j\ge 0$,
\begin{equation}
r_j= \sum_{k=0}^{\lfloor
j/3\rfloor}f_{j-3k,3k}\Big/\sum_{k=0}^{\lfloor
(j-1)/3\rfloor}f_{j-1-3k,3k}.\lab{r6}
\end{equation}
By Corollary~\ref{c:triple}, for all $j=o(n^{1/3})$, $r_j\sim
f_{j,0}/f_{j-1,0}$. By Lemma~\ref{l:single}, this ratio is
asymptotic to $2/j$. This verifies Theorem~\ref{t:Gnm} (a). Let
$\gamma(n)=n/\ln\ln n$. Lemma~\ref{l:largeJ} verifies condition (c). 
The proof is completed by
verifying condition (b). Since $r_j\sim 2/j$ for all $j=o(n^{1/3})$,
we only need to show that for all $n^{1/3}/\ln n\le j\le \gamma(n)$,
$r_j\le m/2N$. It follows directly from the following two facts.

(a) Let $\widehat k=\min\{k: j-3k\le\ln n\}$. By
Corollary~\ref{c:triple},
$$
\sum_{k=0}^{\lfloor j/3\rfloor}f_{j-3k,3k}\sim \sum_{k=0}^{\widehat
k}f_{j-3k,3k},\ \ \ \sum_{k=0}^{\lfloor
(j-1)/3\rfloor}f_{j-1-3k,3k}\sim\sum_{k=0}^{\widehat
k}f_{j-1-3k,3k}.
$$

(b) By Lemma~\ref{l:single}, for all $0\le k\le \widehat k$,
$f_{j-3k,3k}/f_{j-1-3k,3k}=o(1)$. \qed

\section{Proofs of Theorems~\ref{t:Gnm} and~\ref{t:Gnp}}
\lab{sec:proofs}

Before approaching Theorems~\ref{t:Gnm} and~\ref{t:Gnp}, we first
prove a technical lemma.

\begin{lemma}\lab{l:cal} Let $N=\binom{n}{2}$ and let $p=m(n)/N$, where $0<m(n)<N$. Then for any integer $\ell=\ell(n)\ge 0$ such that $\limsup_{n\to \infty}\ell(n)/m(n)<1$,
$$
\binom{N-\ell}{m-\ell}/\binom{N}{m}=p^{\ell}\exp\left(-\frac{1-p}{pN}\frac{\ell^2-\ell}{2}+O(\ell^3/m^2)\right).
$$
Moreover, if $\ell=\Omega(\sqrt{m})$, then
$$
\binom{N-\ell}{m-\ell}/\binom{N}{m}=p^{\ell}\exp\left(-\frac{1-p}{pN}\frac{\ell^2}{2}+O(\ell^3/m^2)\right).
$$
\end{lemma}
\proof
\begin{eqnarray*}
\binom{N-\ell}{m-\ell}/\binom{N}{m}&=&\frac{[m]_{\ell}}{[N]_{\ell}}=\prod_{i=0}^{\ell-1}\frac{m-i}{N-i}\\
&=&\prod_{i=0}^{\ell-1}\frac{m}{N}\exp\left(-\frac{i}{m}+\frac{i}{N}+O(i^2/m^2)\right)\ \ (\mbox{since}\ \limsup_{n\to \infty}\ell(n)/m(n)<1)\\
&=&p^{\ell}\exp\left(-\frac{1-p}{pN}\frac{\ell^2-\ell}{2}+O(\ell^3/m^2)\right).
\end{eqnarray*}
If we have further that $\ell=\Omega(\sqrt{m})$, then
$\ell/pN=O(\ell^3/m^2)$. \qed\ss

\no{\bf Proof of Theorem~\ref{t:Gnm}.
 }  In this proof, the
probability space refers to the random graph $\G(n,m)$ only. Let
$s=|\S|$. By~\eqn{mu-la} and~\eqn{mu},
\begin{eqnarray*}
\ex X_n&=&
s(m/N)^h\exp\left(-\frac{N-m}{mN}\frac{h^2}{2}+O(h^3/m^2)\right).
\end{eqnarray*}
We also have
\begin{eqnarray*}
\ex X_n^2&=&\sum_{j=0}^{h}f_j
\binom{N-(2h-j)}{m-(2h-j)}/\binom{N}{m}.
\end{eqnarray*}
Let $g(j)=f_j \binom{N-(2h-j)}{m-(2h-j)}/\binom{N}{m}$. By condition
(a), for every $K>0$ and any $1\le j\le K h^2/m$,
\begin{equation}
\frac{g(j)}{g(j-1)}=r_j\cdot\frac{N}{m}(1+O(h/m))=\frac{h^2}{mj}(1+O(h/m)+o(m/h^2))=\frac{h^2}{mj}(1+o(m/h^2)),\lab{ratio}
\end{equation}
where the last equality holds because $h^3=o(m^2)$. By condition (c)
and the fact that for any integer $0\le j\le h$,
$\binom{N-(2h-j)}{m-(2h-j)}\le \binom{N-h}{m-h}$, we also have that
$$
\sum_{j>\gamma(n)}g(j)\le
t(n)\binom{N-h}{m-h}/\binom{N}{m}=t(n)\mu_n/s.
$$

 Then   for all sufficiently large $K>0$,
\begin{eqnarray}
\ex
X_n^2&=&\sum_{j=0}^{h}g(j)=\sum_{j=0}^{Kh^2/m}g(j)+O(g(Kh^2/m))+O(t(n)\mu_n/s)\nonumber\\
&=&\left(1+O\left(K^{-1}\right)\right)\sum_{j=0}^{Kh^2/m}g(j)+O(t(n)\mu_n/s),\lab{conv11}
\end{eqnarray}
where the second equality holds because of condition (b) and the
last equality holds by~\eqn{ratio}. 
Next, we estimate $\sum_{j=0}^{Kh^2/m}g(j)$. By~\eqn{ratio} and
Lemma~\ref{l:cal},
\begin{eqnarray}
&&\sum_{j=0}^{Kh^2/m}g(j)=f_0\frac{\binom{N-2h}{m-2h}}{\binom{N}{m}}\sum_{j=0}^{Kh^2/m}\frac{(h^2/m)^j}{j!}(1+o(jm/h^2))\nonumber\\
&&=f_0\cdot
(m/N)^{2h}\exp\left(-\frac{N-m}{mN}\frac{(2h)^2}{2}+O(h^3/m^2)\right)\Big(\exp(h^2/m +o(1))+\Gamma(K)\Big),\nonumber\\
&&= f_0\cdot
(m/N)^{2h}\exp\left(-\frac{N-m}{mN}2h^2\right)\exp(h^2/m)\Big(1+o(1)+O(\Gamma(K)\exp(-h^2/m))\Big),
\lab{conv22}
\end{eqnarray}
where
$$
\Gamma(K)=O\left(\frac{(h^2/m)^{Kh^2/m}}{(Kh^2/m)!}\right)=O\left(\left(\frac{(eh^2/m)}{(Kh^2/m)}\right)^{Kh^2/m}\right),
$$
which goes to $0$ as $K\to\infty$, since $h^2/m=\Omega(1)$. By~\eqn{conv11} and~\eqn{conv22}, for every $\eps>0$, there is a sufficiently large $K$, such that
 \begin{equation}
\ex X_n^2=(1+O(\eps)) f_0\cdot
(m/N)^{2h}\exp\left(-\frac{N-m}{mN}2h^2\right)\exp(h^2/m)+O(t(n)\mu_n/s).\lab{sqrX}
\end{equation}
We also have
$$
s^2=\sum_{j=0}^{h}f_j=f_0\sum_{j=0}^h\prod_{i=1}^jr_i.
$$
With the same reasoning as before, it is enough to sum over the
first $Kh^2/N$ terms, leaving an arbitrarily small tail plus an error term
$O(t(n))$. This yields
$$
 s^2=(1+O(\eps))
f_0\exp(h^2/N)+O(t(n)).
$$
Since $t(n)=o(\mu_ns)=o(s^2)$ by condition (c), we obtain
$$
f_0=(1+O(\eps)) s^2\exp(-h^2/N).
$$
Combining with~\eqn{sqrX} and again by condition (c), we obtain
\begin{eqnarray*}
\ex X_n^2&=&(1+O(\eps))
s^2(m/N)^{2h}\exp\left(-\frac{N-m}{mN}2h^2\right)\exp(h^2/m-h^2/N)+O(t(n)\mu_n/s)\\
&=&(1+O(\eps))s^2(m/N)^{2h}\exp\left(-\frac{N-m}{mN}h^2\right)+o(\mu_n^2)=
(1+O(\eps))(\ex X_n)^2.
\end{eqnarray*}
As this holds for every $\eps>0$, we have $\ex X_n^2=(1+o(1))(\ex X_n)^2$.
Then for
every $\eps>0$,
$$
\pr(|X_n/\ex X_n-1|>\eps)\to 0, \ \ \mbox{as}\ n\to\infty,
$$
by Chebyshev's inequality. \qed \ss

 \no{\bf Proof of
Theorem~\ref{t:Gnp}.\ } Let $Y_n$ denote the number of edges in
$\G(n,p)$, then $Y_n\sim Bin(N,p)$. Hence we have
\begin{equation}
Y_n-pN=O_p(\sqrt{p(1-p)N}),\lab{concentration}
\end{equation}
where $f(n)=O_p(g(n))$ for some $g(n)\ge 0$ means
$\pr(|f(n)|>Kg(n))\to 0$ as $K\to\infty$ and $n\to\infty$. Similarly
we use the notation $f(n)=o_p(g(n))$ meaning that for every
$\eps>0$, $\pr(|f(n)|>\eps g(n))\to 0$ as $n\to \infty$. Since
$X_n/\ex_{\G(n,m)X_n}\xrightarrow{p} 1$ in $\G(n,m)$ for all
$m=pN+O(\sqrt{p(1-p)N})$ by assumption and
$\ln(\ex_{\G(n,m)}X_n)=\ln|\S|+h\ln(m/N)+(N-m)h^2/2mN+o(1)$
by~\eqn{mu}, \remove{

*********
Assume further that $s$ is sufficiently large so that
condition~\eqn{cond1} holds. Then for all integers $m$ such that
$m=pN+O(\sqrt{p(1-p)N})$, all assumptions in Theorem~\ref{t:Gnm}
satisfy and $\ex_{\G(n,m)} X_n\to\infty$ as $n\to\infty$. *********}
  by conditioning on $Y_n$, we have
\begin{equation}
\ln X_n-\ln
|\S|-h\ln(Y_n/N)+\frac{1-Y_n/N}{Y_n}\frac{h^2}{2}\xrightarrow{p}
0.\lab{conv1}
\end{equation}
By~\eqn{concentration},
\begin{eqnarray}
\frac{1-Y_n/N}{Y_n}\frac{h^2}{2}&=&\frac{h^2(1-p)}{2Np}\left(1+O_p\left(\sqrt{\frac{p}{(1-p)N}}+\sqrt{\frac{1-p}{pN}}\right)\right)=\frac{h^2(1-p)}{2Np}+o_p(1),\lab{conv2}
\end{eqnarray}
where the equality above holds because $h^3=o(p^2n^4)$. We also have
\begin{equation}
\ln(Y_n/N)=\ln p(1+Y^*_n\sqrt{(1-p)/pN})=\ln
p+\sqrt{(1-p)/pN}Y^*_n+O_p((1-p)/pN),\lab{conv3}
\end{equation}
where
$$
Y^*_n=\frac{Y_n-pN}{\sqrt{p(1-p)N}}
$$
is the normalised variable of $Y_n$. Recall that $\la_n=|\S|p^h$
from~\eqn{mu-la} and $\ex X_n=\la_n$. Combining
with~\eqn{conv1}--\eqn{conv3}, we have
\begin{equation}
\ln(X_n/\la_n)+\frac{\beta_n^2}{2}=\beta_n Y^*_n+o_p(1).\lab{conv4}
\end{equation}
Since $\beta_n=\Omega(1)$,~\eqn{conv4} immediately yields
$$
\frac{\ln(e^{\beta_n^2/2}X_n/\la_n)}{\beta_n}= Y^*_n + o_p(1).
$$
Since $Y^*_n\xrightarrow{d} \N(0,1)$, the theorem follows.\qed \ss

\section{Concluding remarks}
\lab{sec:concluding}

 It was proved in~\cite{J3} that $m>>n^{3/2}$ is required for the concentration of $X_n$ in $\G(n,m)$, where $X_n$ denotes
 the number of Hamilton cycles or perfect matchings or spanning trees, as the variable will become asymptotically log-normally
 distributed when $m=\Theta(n^{3/2})$. However, we do not think this condition is sufficient in the case of triangle-factors. It is surprising that the critical
  point of $m$ when $X_n$ changes from small deviation ($\ex X_n^2\sim (\ex X_n)^2$) to large deviation
  ($\limsup_{n\to\infty} \ex X_n^2/(\ex X_n)^2>1$) in $\G(n,m)$ seems to be different for Hamilton cycles and
   for triangle-factors. We guess $m=n^{5/3}$ might be the critical point for the latter case.

As explained in Section~\ref{sec:hamilton}, the most interesting set
$\S$ to be studied is perhaps the one containing graphs isomorphic
to an unlabelled graph $H_n$ on $n$ vertices. Unfortunately, it is not easy to
define the sequence $(H_n)_{n\ge 1}$ in general and for a
general $H_n$, computing $r_j$ might be hard. It will be
interesting to discover more classes of such graph sequences $(H_n)$ and see whether the
corresponding random variables $X_n$ follow the log-normal paradigm.

\end{document}